\documentclass[11pt]{article}
\usepackage{amsmath}
\usepackage{amsthm}
\usepackage{amssymb}

\newtheorem{theo}{Theorem}[section]
\newtheorem{definition}{Definition}[section]
\newtheorem{prop}[theo]{Proposition}
\newtheorem{lemma}[theo]{Lemma}

\newcommand{\FF}{{\cal F}}

\begin{document}
\date{}

\title{
Clique coloring of dense random graphs
}

\author{Noga Alon
\thanks{Sackler School of Mathematical Sciences
and Blavatnik School of
Computer Science, Tel Aviv University, Tel Aviv 6997801, Israel,
and CMSA, Harvard University, Cambridge, MA 02138, USA.
Email: {\tt nogaa@tau.ac.il}.  Research supported in part by
an ISF grant and by a GIF grant.
}
\and
Michael Krivelevich
\thanks{
Sackler School of Mathematical Sciences,
Tel Aviv University, Tel Aviv 6997801, Israel.
Email: {\tt krivelev@post.tau.ac.il}.
Research supported in part by USA-Israeli BSF grant 2014361
and by ISF grant 1261/17.}
}

\maketitle
\begin{abstract}
The clique chromatic number of a
graph $G=(V,E)$ is the minimum number of colors in a
vertex coloring so that no maximal (with respect to containment)
clique is monochromatic.
We prove that the
clique chromatic
number of the binomial random graph $G=G(n,1/2)$ is, with high
probability, $\Omega( \log n)$. This settles a problem of
McDiarmid, Mitsche and Pra{\l}at who proved that it is
$O( \log n)$ with high probability.
\end{abstract}

\section{The main result}

A {\em clique} in an undirected graph $G=(V,E)$ is {\em maximal}
if it is not
properly contained in a larger clique. A {\em clique coloring} of $G$ is
a vertex coloring so that no maximal clique (with at least two
vertices) is monochromatic. Let
$\chi_c(G)$ denote the minimum possible number of colors in
a clique coloring of $G$. This invariant is called the {\em clique
chromatic number} of $G$ and has been studied in a considerable
number of papers, see
\cite{AST,BGGPS,CDM,CK,De,GHM,KM,KT,LSK,MMP,MS,SLK}. McDiarmid, Mitsche and Pra{\l}at \cite{MMP} initiated the study
of $\chi_c(G)$ for the random binomial graph $G=G(n,p)$.
While for sparse random graphs their upper and lower bounds for the
typical behavior of $\chi_c(G(n,p))$ are rather close to each
other, for the dense case they do not have any nontrivial lower bound.
In particular, for the random graph $G=G(n,1/2)$ they proved that
with high probability (whp, for short), that is, with probability
tending to $1$ as $n$ tends to infinity, $\chi_c(G) \leq
(1/2+o(1)) \log_2 n$ and raised the problem of proving any
nontrivial lower bound. In this note we show that the logarithmic
estimate is tight, up to a constant factor.
\begin{theo}
\label{t11}
There exists an absolute positive constant
$c$ so that whp the random graph $G=G(n,1/2)$ satisfies
$\chi_c(G) \geq c \log_2 n$. Therefore, whp
$\chi_c(G)=\Theta(\log n)$.
\end{theo}
The proof appears in the next two sections. Throughout the proof we
assume, whenever this is needed, that $n$ is sufficiently large. All
logarithms are in base $2$, unless otherwise  specified. To
simplify the presentation, we omit all floor and ceiling signs
whenever these are not crucial, and make no attempt to optimize
the absolute constants in our estimates.

\section{Preliminaries}
This section includes the main technical part of the proof. It introduces several typical properties of the random graph $G(n,1/2)$,
established in the following three lemmas. The
proofs of the first two are straightforward, while that of the third one
requires some work.
\begin{lemma}
\label{l21}
Let $G=G(n,1/2)=(V,E)$.
Then whp for every set $S$ of at most $\frac{1}{2000}
\log n$ vertices of $G$ there are more than
$n^{0.999} \log n $ vertices in $V-S$ that are are not adjacent to
any vertex of $S$.
\end{lemma}
\begin{lemma}
\label{l22}
The following holds for the random graph $G=G(n,1/2)$ whp.
For every set $Y$ of $|Y|=y \geq n^{0.999}$ vertices of $G$,
the number of vertices in $V-Y$ that have at most
$0.41y$ non-neighbors in $Y$ is smaller than $\frac{1}{4} \log n$.
\end{lemma}
The third lemma is a bit more technical. It is convenient to define first
the following property of a set of vertices $Y$.
\begin{definition}
\label{d24}
Let $Y \subsetneq V$ be a set of vertices of an $n$-vertex graph $G=(V,E)$.
Call $Y$ {\em significant} if
it satisfies the following two conditions:
\begin{itemize}
\item
Every vertex $v \in V-Y$ has at least $n^{0.999}$ non-neighbors in
$Y$;
\item
The number of vertices $v \in V-Y$
that have at most $0.41y$ non-neighbors in
$Y$ is at most $\frac{1}{4} \log n$.
\end{itemize}
\end{definition}
Note that, by definition, every significant set has size at least
$n^{0.999}$.
\begin{lemma}
\label{l23}
Let $G=G(n,1/2)=(V,E)$. Then whp every significant set
$Y$ in $G$ contains a clique $K$ of size $k=1.9 \log n$
so that every vertex $v \in V-Y$ has at least one non-neighbor
in $K$.
\end{lemma}

We proceed with the proofs of the three lemmas above.
The proofs of the first two are very easy.
\vspace{0.2cm}

\noindent
{\bf Proof of Lemma \ref{l21}:}\,
Fix a set $S$  of $s \leq \frac{1}{2000} \log n$
vertices. The number of vertices in $V-S$ that are not adjacent to
any member of $S$ is a binomial random variable with parameters
$n-s$ and $1/2^{s}$ whose expectation is $(n-s)2^{-s}
\geq (1-o(1)) n^{0.9995}$. By the standard estimates for binomial
distributions (c.f., e.g., \cite{AS}, Theorem A.1.13) the
probability that this number is smaller than $n^{0.999} \log n$,
which is less than half its expectation, is smaller than
$e^{- n^{0.9995} /9}$. Note that the number of possible sets $S$ is
(much) less
than $n^{\log n}$, since
$$
\sum_{s=0}^{\log n/2000} {n \choose s} \leq
\sum_{s=0}^{\log n/2000} n^s \leq  2n^{\log n/2000} < n^{\log n}.
$$
The desired result thus follows, by the union bound.
\hfill $\Box$
\vspace{0.2cm}

\noindent
{\bf Proof of Lemma \ref{l22}:}\,
If $G$ contains a set $Y$ of size $|Y|=y\ge n^{0.999}$ violating
the lemma's claim, then there is a subset $X\subset V\setminus Y$ of
size $|X|=x=\frac{1}{4}\log n$ such that every vertex $x\in X$ sends
at least $0.59y$ edges to $Y$. This implies that $G$ has at least
$0.59xy$ edges crossing between $X$ and $Y$. For two given sets $X,Y$
as above, the number of edges $e(X,Y)$ between $X$ and $Y$ in $G(n,1/2)$
is distributed binomially with parameters $xy$ and $1/2$. Using, again,
the known estimates for binomial distributions (c.f. \cite{AS}, Theorem
A.1.1), we obtain that the probability that $e(X,Y)\ge 0.59xy$ is at
most $e^{-2 \cdot 0.09^2 xy}<e^{-0.016xy}$. Summing over all possible
choices of $X$ and $Y$, it follows that the probability of the existence
of a set $Y$ violating the lemma's claim is at most
\begin{eqnarray*}
&&\sum_{y\ge n^{0.999}}\binom{n}{y}\binom{n-y}{x}e^{-0.016xy}
\le \sum_{y\ge n^{0.999}} \left(\frac{en}{y}\right)^yn^xe^{-0.016xy}\\
&\le&  \sum_{y\ge n^{0.999}} n^x\cdot \left(en^{0.001}\cdot
e^{-\frac{0.016\log n}{4}}\right)^y
\le \sum_{y\ge n^{0.999}} n^{\log n/4}
\cdot n^{-0.002\cdot n^{0.999}}= o(1)\,,
\end{eqnarray*}
completing the proof of the lemma. \hfill $\Box$
\vspace{0.2cm}

\noindent
{\bf Proof of Lemma \ref{l23}:}\,
Fix a set $Y$ of $y \geq n^{0.999}$ vertices of $G=G(n,1/2)=(V,E)$
and expose all edges of $G$ between $Y$ and $V-Y$. If $Y$ is not
significant then there is nothing to prove, we thus assume that
$Y$ is significant. Put $r=\frac{1}{4} \log n$ and $s=k-r=
1.65 \log n$. Let $B$ be a set of exactly
$r$ vertices in $V-Y$ containing
all vertices in $V-Y$ that have
less than $0.41y$ non-neighbors in $Y$. Put $m=n^{0.99}$ and
choose, for each $v \in B$, a subset of $Y$ of size $m$ consisting
of non-neighbors of $v$, where all these subsets are pairwise
disjoint. (It is easy to choose these sets sequentially, as
each $v \in B$ has at least $n^{0.999}$ non-neighbors in $Y$.)
This defines $r$ subsets which we denote by
$Z_1,Z_2, \ldots Z_r$. Put
$Y'=Y - \cup_{i=1}^r  Z_i$, $|Y'|=y'$, and note that each vertex
$v \in V-(Y \cup B)$ has at least $0.41y -r n^{0.99}> 0.405y'$
non-neighbors in $Y'$.
\vspace{0.1cm}

\noindent
{\bf Claim: }\,
There are $s$ pairwise disjoint subsets $Z_{r+1}, Z_{r+2}, \ldots
,Z_{k}$ of $Y'$, each of size exactly $m$, so that every vertex
$v \in V-(Y \cup B)$ has at least $0.4m$ non-neighbors in each of
the subsets $Z_j,~r+1 \leq j \leq k$.
\vspace{0.1cm}

\noindent
{\bf Proof of claim:}\, Choose the sets randomly  and apply the
standard estimates for hypergeometric distributions
(c.f. \cite{JLR}, Theorem
2.10). \hfill $\Box$
\vspace{0.1cm}

\noindent
Let $\FF$ be the family of all  subsets of size $k$ of
$Y$ that contain exactly one element in each set $Z_i$  and
contain at least one non-neighbor  of each vertex $v \in V-Y$.
Note that by the definition of the first sets $Z_1,\ldots
Z_r$, each set that contains an element from each $Z_i$ has at least
one non-neighbor of  each vertex $v \in B$. On the other hand, for
each fixed $v \in V-(Y \cup B)$, when we choose randomly
one member from each $Z_j$ for $r+1 \leq j \leq k$, the probability
that we do not choose any non-neighbor of $v$ is at most
$0.6^{s}=0.6^{1.65 \log n} <\frac{1}{n^{1.1}}$. Therefore, by the
union bound, almost all of these choices do include at least one
non-neighbor of each such $v$ and hence
$$
|\FF| \geq (1-o(1)) m^k =(1-o(1)) n^{0.99 \cdot 1.9 \log n}.
$$

We now expose the edges in the induced subgraph of $G$ on
$\cup _j Z_j$ and show that the probability that none of the members of $\FF$ is a clique is much smaller than $2^{-n}$. This can be proved in several ways, either by using martingales (see
\cite{AACKR}, Section 4.1 for a similar argument), or by using
Talagrand's Inequality, or by using the extended
Janson's Inequality (c.f., \cite{AS}, Theorem 8.1.2). The last
alternative seems to be the shortest, and we proceed with its
detailed description.

For each member $K$ of $\FF$, let $x_K$ denote the indicator random
variable whose value is $1$ if and only if $K$ is a clique in $G$
and let $X=\sum_{K \in \FF} X_K$. Our objective is to show that
$X>0$ with probability that is close enough to 1 to enable applying the union bound over all relevant sets $Y$.
The expectation of each $X_K$ is clearly
$$
E(X_K)=2^{-{k \choose 2}}.
$$
Thus, by linearity of expectation,
$$
E(X) =|\FF|2^{-{k \choose 2}} =(1-o(1)) m^k 2^{-{k \choose 2}}
=(1-o(1))(m2^{-(k-1)/2})^k > n^{0.03 k} >n^{10}
$$
with (a lot of) room to spare.

Put $\mu=E(X)$, and define
$\Delta = \sum_{K,K'} \mbox{Prob}[X_K=X_{K'}=1]$
where the summation is over all
(ordered) pairs $K,K'$ of members of $\FF$  that satisfy
$2 \leq |K \cap K'| \leq k-1$.
By the extended Janson Inequality the probability that $X=0$ is at
most $e^{-\mu^2/2\Delta}.$

Note that $\Delta=\sum_{i=2}^{k-1} \Delta_i$ where $\Delta_i$ is
the contribution of pairs $K,K'$ with $K,K' \in \FF$, $|K \cap
K'|=i$. Thus
$$
\Delta_i \leq |\FF| 2^{-{k \choose 2}} {k \choose i} (m-1)^{k-i}
2^{-{k \choose 2}+{i \choose 2}} \leq
m^k 2^{-2{k \choose 2}+{i \choose 2}} {k \choose i} m^{k-i}.
$$
We next prove that
\begin{equation}
\label{e21}
\Delta=\sum_{i=2}^{k-1}\Delta_i \leq (1+o(1))\frac{k^2}{m^2} \mu^2.
\end{equation}
To do so, consider the following cases.
\vspace{0.1cm}

\noindent
{\bf Case 1:}\,  $i=2$. In this case
$$
\frac{\Delta_2}{\mu^2} \leq (1+o(1)) \frac{{k \choose 2}\cdot 2}{m^2}
\leq (1+o(1)) \frac{k^2}{m^2}.
$$
\vspace{0.1cm}

\noindent
{\bf Case 2:}\,  $3  \leq i < 100$. Here
$$
\frac{\Delta_i}{\mu^2} \leq (1+o(1)) \frac{{k \choose i} 2^{{i
\choose 2}}}{m^i} <(1+o(1))\left(\frac{k2^{i/2}}{m}\right)^i \leq
\left(\frac{k2^{50}}{m}\right)^3 =\frac{1}{m^{3-o(1)}}.
$$
\vspace{0.1cm}

\noindent
{\bf Case 3:}\,  $100 \leq i \leq k-2$. In this case
$$
\frac{\Delta_i}{\mu^2} \leq (1+o(1)) \frac{{k \choose i} 2^{{i
\choose 2}}}{m^i} <(1+o(1))\left(\frac{k2^{i/2}}{m}\right)^i \leq
\left(\frac{1}{n^{0.04-o(1)}}\right)^{100} =\frac{1}{n^{4-o(1)}}.
$$
Summing the contributions for all $i$, $2 \leq i \leq k-1~(< \log
n)$, the inequality (\ref{e21}) follows.

By the extended Janson Inequality this implies that the probability
that $X=0$ is at most
$$
e^{-\mu^2/2\Delta} \leq e^{-(1+o(1))m^2/2k^2}<e^{-n^{1.98-o(1)}}.
$$
As the number of possible significant sets $Y$ is smaller than $2^n$,
the assertion of the lemma follows, by the union bound. \hfill
$\Box$

\section{Completing the proof}

In this section we prove Theorem \ref{t11}. By the results in the
previous section it suffices to prove the following deterministic
statement.
\begin{prop}
\label{p31}
Let $G=(V,E)$ be a graph on $n$ vertices satisfying the assertions of
Lemma \ref{l21}, Lemma \ref{l22} and Lemma \ref{l23}. Then
$\chi_c(G) > \frac{1}{2000} \log n$.
\end{prop}

\noindent
{\bf Proof:}\, Assume this is false, and let $Y_1,Y_2, \ldots ,Y_s$
be a partition of the vertex set $V$ into disjoint non-empty sets, each
containing no maximal clique of $G$, where $s  \leq \frac{1}{2000}
\log n$. It is easy to verify that the conclusions
of Lemmas \ref{l21}, \ref{l22} guarantee in particular
that $G$ contains at least one edge, and thus $s>1$.
For each $i$, $1 \leq i \leq s$, let $v_i \in V-Y_i$ be a
vertex with the minimum number of non-neighbors in $Y_i$ among all
vertices in $V-Y_i$, and let $t_i$ denote the number of these
non-neighbors. By definition, the number of vertices of $G$ which are
not adjacent to any vertex of $S=\{v_1, v_2, \ldots ,v_s\}$ is
at most  $\sum_{i=1}^s t_i$, and since
$G$ satisfies the assertion of Lemma \ref{l21} this number
exceeds $n^{0.999} \log n$. By averaging there exists an index $i$
so that $t_i \geq n^{0.999}$. Fix such an $i$ and note that
$Y=Y_i$ is a significant set. Indeed, by the
definition of $t_i$ each $v \in V-Y_i$  has at least
$t_i \geq n^{0.999}$ non-neighbors in $Y_i$, and as $G$ satisfies
the conclusion of Lemma \ref{l22} the number of vertices
$v \in V-Y_i$ that have at most $0.41|Y_i|$ non-neighbors in $Y_i$
is at most $\frac{1}{4} \log n$.

Since $G$ satisfies the conclusion of Lemma \ref{l23}, the set  $Y:=Y_i$
contains a clique $K$ of size $1.9 \log n$ which contains at least one
non-neighbor of each vertex $v \in V-Y_i$. This clique is contained
in a maximal clique of $G$, call it $K'$. However, $K'$ is
contained in $Y_i$, as $K$ has a non-neighbor of each vertex $v
\in V-Y_i$. Thus $K'$ is a maximal clique of $G$ which is
contained in $Y_i$, contradiction. This completes the proof of the
proposition and hence of Theorem \ref{t11}. \hfill $\Box$

\section{Concluding remarks and open problems}

We have shown that the clique chromatic number $\chi_c(G)$ of the random
graph $G=G(n,1/2)$ is, whp, $\Theta(\log n)$.  The same proof applies
to binomial random graphs with any constant edge probability bounded away
from $0$ and $1$. Together with the upper bound proved in \cite{MMP} we
conclude that for any fixed $p$, $0 < p <1$, the random graph $G=G(n,p)$
satisfies, whp, $\chi_c(G)=\Theta_p(\log n)$. Note, however, that as
$p$ tends to $0$ the ratio between the upper bound of
\cite{MMP}, which is $\Theta(\log n/p)$, and our lower bound increases.
Indeed, the natural modification of the parameters in our proof here for
small values of $p$ bounded away from zero provides only an
$\Omega(\log n/ \log (1/p))$ lower bound.

We have made no attempt to optimize the absolute constants in
our estimates. The constant $\frac{1}{2000}$
can certainly be improved, but
it seems that the method as it is does not suffice to determine
the tight constant here. It seems plausible that
$$
\chi_c(G(n,1/2))=(1/2+o(1))\log n
$$
whp.

\section*{Acknowledgments}
We thank Colin McDiarmid for helpful discussions, Andrzej Dudek, Tomasz
\L uczak and Pawe\l\ Pra\l at for their input, and two anonymous
referees for their comments. This project was
initiated during a workshop in
Monash University, Melbourne, in June 2016. We thank the organizers
of this workshop for their
hospitality.


\begin{thebibliography}{99}
\bibitem{AACKR}
N. Alon, V. Asodi, C. Cantor, S. Kasif and J. Rachlin,
Multi-node graphs: a framework for multiplexed biological assays,
Journal of Computational Biology 13, 2006, 1659--1672.
\bibitem{AS}
N. Alon and J. H. Spencer,  {\bf The Probabilistic
Method}, Fourth Edition, Wiley, 2016.
\bibitem{AST}
T. Andreae, M. Schughart and Zs. Tuza, Clique-transversal sets of
line graphs and complements of
line graphs, Discrete Math. 88, 1991, 11--20.
\bibitem{BGGPS}
G. Bacs\'o, S. Gravier, A. Gy\'arf\'as, M. Preissmann and A. Seb\'o,
Coloring the maximal cliques of graphs,
SIAM Journal on Discrete Mathematics 17, 2004, 361--376.
\bibitem{CDM}
C.N. Campos, S. Dantas and C.P. de Mello, Colouring
clique-hypergraphs of circulant graphs, Electron.
Notes Discret. Math. 30, 2008, 189--194.
\bibitem{CK}
M.R. Cerioli and A.L. Korenchendler, Clique-coloring circular-arc
graphs, Electron. Notes Discret.
Math. 35, 2009, 287--292.
\bibitem{De}
D. D\'efossez, Clique-coloring some classes of odd-hole-free
graphs, J. Graph Theory 53, 2006,
233--249.
\bibitem{GHM}
S. Gravier, C. Ho\'ang and F. Maffray, Coloring the hypergraph of
maximal cliques of a graph with no
long path, Discrete Math., 272, 2003, 285--290.
\bibitem{JLR}
S. Janson, T. \L uczak and A. Ruci\'nski, {\bf Random graphs},  Wiley, 2000.
\bibitem{KM}
S. Klein and A. Morgana, On clique-colouring with few P4s, J.
Braz. Comput. Soc. 18, 2012, 113--119.
\bibitem{KT}
J. Kratochv\'il and Zs. Tuza, On the complexity of bicoloring
clique hypergraphs of graphs, J. Algorithms,
45, 2002, 40--54.
\bibitem{LSK}
Z. Liang, E. Shan and L. Kang, Clique-coloring claw-free graphs,
Graphs and Combinatorics,
2015, 1--16.
\bibitem{MMP}
C. McDiarmid, D. Mitsche and P. Pra{\l}at,
Clique coloring of binomial random graphs, arXiv 1611.01782
\bibitem{MS}
B. Mohar and R. Skrekovski,
The Gr\"otzsch Theorem for the
hypergraph of maximal cliques, Electronic
J. Comb. 6, 1999, R26.
\bibitem{SLK}
E. Shan, Z. Liang and L. Kang, Clique-transversal sets and
clique-coloring in planar graphs, European
J. Combin. 36, 2014, 367--376.
\end{thebibliography}
\end{document}